\documentclass[reqno,11pt]{amsart}

\usepackage{xcolor}
\usepackage{graphicx}
\usepackage[hidelinks]{hyperref}
\usepackage{amscd}

\usepackage{amsmath}
\usepackage{amsthm}
\usepackage{amssymb}
\usepackage{latexsym,array}
\usepackage{amsfonts}
\usepackage{shadow}
\usepackage{amsbsy}
\usepackage{amssymb}
\usepackage{dsfont}
\usepackage{mathtools}
\usepackage{enumitem}
\usepackage{doi}
\usepackage{color}
\usepackage{graphicx}
\usepackage[hidelinks]{hyperref}
\usepackage{cleveref}
\usepackage{fullpage}
\usepackage{comment}

\usepackage{dsfont}

\usepackage{cleveref}

\numberwithin{equation}{section}

\newcommand{\PRes}{\mathcal{Q}}

\newtheorem{theorem}{Theorem}[section]

\theoremstyle{definition}
\newtheorem{remark}[theorem]{{\bf Remark}}
\newtheorem{definition}[theorem]{Definition}

\newcommand{\cc}{\mathbb{C}}

\newcommand{\rr}{\mathbb{R}}

\renewcommand{\Re}{\mathrm{Re}}

\usepackage{amscd}
\usepackage{tikz-cd}
\usepackage{tikz}

\crefname{enumi}{}{}
\crefname{enumii}{}{}

\title[]{An introduction to the fine structures on the $S$-spectrum}
\author[F. Colombo]{Fabrizio Colombo}
\address{(FC)
	Politecnico di Milano\\Dipartimento di Matematica\\Via E. Bonardi, 9\\20133
	Milano, Italy}
\email{fabrizio.colombo@polimi.it}

\author[A. De Martino]{Antonino De Martino}
\address{(ADM)
	Politecnico di Milano\\Dipartimento di Matematica\\Via E. Bonardi, 9\\20133
	Milano, Italy
} \email{antonino.demartino@polimi.it}

\author[S. Pinton]{Stefano Pinton}
\address{(SP)
	Politecnico di Milano\\Dipartimento di Matematica\\Via E. Bonardi, 9\\20133
	Milano, Italy
} \email{stefano.pinton@polimi.it}

\author[I. Sabadini]{Irene Sabadini}
\address{(IS)
	Politecnico di Milano\\Dipartimento di Matematica\\Via E. Bonardi, 9\\20133
	Milano, Italy
} \email{irene.sabadini@polimi.it}

\author[P. Schlosser]{Peter Schlosser}
\address{(PS)
	Politecnico di Milano\\Dipartimento di Matematica\\Via E. Bonardi, 9\\20133
	Milano, Italy
} \email{pschlosser@math.tugraz.at}

\date{}

\begin{document}

	\maketitle

\begin{abstract}

Holomorphic functions are fundamental in operator theory and their Cauchy formula is a crucial tool for defining functions of operators.
The Fueter-Sce extension theorem (often called Fueter-Sce mapping theorem)  provides a two-step procedure for extending holomorphic functions to hyperholomorphic functions. In the first step, slice hyperholomorphic functions are obtained, and their associated Cauchy formula establishes the $S$-functional calculus for noncommuting operators on the $S$-spectrum. The second step produces axially monogenic functions, which lead to the development of the monogenic functional calculus.
In this review paper we discuss the second operator in the Fueter-Sce mapping theorem that takes
slice hyperholomorphic to axially monogenic functions. This operator admits several factorizations which  generate various function spaces and their corresponding functional calculi, thereby forming the so-called fine structures of spectral theories on the $S$-spectrum.

\end{abstract}

\vspace*{+4mm}
\centerline{\bf This paper is dedicated to Paula Cerejeiras on the occasion of her 60th Birthday.}

\section{Introduction}\label{INTRO}
Driven by motivations in quaternionic quantum mechanics (see \cite{BF}),
since 2006, the spectral theory of the $S$-spectrum for quaternionic operators
has significantly advanced, with key references including \cite{6COFBook,6ACSBOOK,FJBOOK,CGKBOOK,6JONAME}.
Recent progress within the Clifford algebra framework \cite{6css} has led to the proof of the spectral theorem for fully Clifford operators \cite{CLIFST}.
Furthermore, the $S$-functional calculus
has been successfully extended beyond the confines of Clifford algebras \cite{UNIV}. For the latest developments in interpolation theory, see \cite{COLSCH},
and for recent applications to vector operators in a Clifford algebra, see \cite{CMS24}.

In this note, we aim to provide an accessible introduction to a new direction in the spectral theory related to the $S$-spectrum, specifically focusing on the so-called fine structures. These fine structures arise from the factorization of the second map $T_{FS2}$ in the Fueter-Sce mapping theorem, an important technique for extending holomorphic functions to hyperholomorphic functions, which also includes vector fields as a special case.
The integral representation formulas for these functions, induced by the factorization of $T_{FS2}$, serve as the foundation for defining new functional calculi for both bounded and unbounded operators.

To contextualize our results, we briefly outline the Fueter-Sce mapping theorem, which involves a two-step procedure for extending holomorphic functions of one complex variable to the hypercomplex setting. The first step generates slice hyperholomorphic functions, while the second step produces monogenic functions that lie in the kernel of the Dirac operator.

To formalize this in Section \ref{SEC2} we fix some notations. Let $\mathbb{R}_n$ denote the real Clifford algebra over $n$ imaginary units $e_1, \ldots, e_n$, satisfying the relations  $e_\ell e_m + e_m e_\ell = 0$ for $\ell \neq m$ and $e_\ell^2 = -1$. An element in the Clifford algebra is represented as $\sum_A e_A x_A$, where $A = \{ \ell_1, \ldots, \ell_r \} \in \mathcal{P}\{1, 2, \ldots, n\}$ with $\ell_1 < \cdots < \ell_r$, and $e_A = e_{\ell_1} e_{\ell_2} \cdots e_{\ell_r}$ with $e_\emptyset =e_0= 1$.
A point $(x_0, x_1, \ldots, x_n) \in \mathbb{R}^{n+1}$ is identified with the element
$x = x_0 + \underline{x} = x_0 + \sum_{j=1}^n x_j e_j \in \mathbb{R}_n,$
referred to as a paravector. The real part of $x$, denoted by $\Re(x)$, is $x_0$, and the vector part is $\underline{x} = x_1 e_1 + \cdots + x_n e_n$. The conjugate of $x$ is $\overline{x} = x_0 - \underline{x}$, and the Euclidean modulus of $x$ is given by $|x| = (x_0^2 + \cdots + x_n^2)^{1/2}$. The sphere of purely imaginary vectors with modulus 1 is defined by
$\mathbb{S} := \{\underline{x} = e_1 x_1 + \cdots + e_n x_n \mid x_1^2 + \cdots + x_n^2 = 1\}.$
If $I \in \mathbb{S}$, then $I^2 = -1$, making $I$ an imaginary unit. We denote by
$\mathbb{C}_I = \{u + Iv \mid u, v \in \mathbb{R}\},$
an isomorphic copy of the complex numbers. Given an element $x\in\rr^{n+1}$, we define
\begin{align*}
	[x]=\{y\in\rr^{n+1}\ :\ y={\rm Re}(x)+I |\underline{x}|,\, I\in \mathbb{S}\}.
\end{align*}
A set $U \subseteq \mathbb{R}^{n+1}$ is called {\em axially symmetric} if $[x]\subset U$  for any $x \in U$.
Similar notations apply to the real algebra $ \mathbb{R}_2$, which is known as the algebra of quaternions and is also denoted by $ \mathbb{H}$.
\begin{definition}[Slice Cauchy domain]
	An axially symmetric open set $U\subset \mathbb{R}^{n+1}$ is called a slice Cauchy domain,
	if $U\cap\cc_I$ is a Cauchy domain in $\cc_I$ for any $I\in\mathbb{S}$.
	More precisely, $U$ is a slice Cauchy domain if for any $I\in\mathbb{S}$
	the boundary ${\partial( U\cap\cc_I)}$ of $U\cap\cc_I$ is the union a finite number of non-intersecting piecewise continuously differentiable Jordan curves in $\cc_{I}$.
\end{definition}
\begin{definition}[Slice hyperholomorphic functions (or slice monogenic functions)]\label{sh}
	Let $U\subseteq \mathbb{R}^{n+1}$ be an axially symmetric open set and let
	$$
	\mathcal{U} = \{ (u,v)\in\rr^2: u+ \mathbb{S} v\subset U\}.
	$$
	A function $f:U\to \mathbb{R}^{n+1}$ is called a left
	slice function (or of left axial type), if it is of the form
	\begin{equation}
		\label{form}
		f(x) = \alpha(u,v) + I\beta(u,v)\qquad \text{for } x = u + I v\in U
	\end{equation}
	where the functions $\alpha, \beta: \mathcal{U}\to \mathbb{R}_{n}$ that satisfy the compatibility condition
	\begin{equation}
		\label{EO}
		\alpha(u,v)=\alpha(u,-v), \qquad \beta(u,v)=- \beta(u,-v), \qquad \forall (u,v) \in \mathcal{U}
	\end{equation}
	If in addition $\alpha$ and $\beta$ satisfy the Cauchy-Riemann-equations
	\begin{equation}
		\label{CR}
		\partial_u \alpha(u,v)- \partial_v \beta(u,v)=0, \quad \partial_v \alpha(u,v)+ \partial_u \beta(u,v)=0,
	\end{equation}
	then $f$ is called left slice hyperholomorphic or left slice monogenic.
	A function $f:U\to \mathbb{R}_n$ is called a right slice function if it is of the form
	\begin{align*}
		f(x) = \alpha(u,v) + \beta(u,v)I\qquad \text{for } x = u + I v\in U
	\end{align*}
	where $ \alpha$ and $ \beta: \mathcal{U}\to \mathbb{R}_{n}$ satisfy \eqref{EO}.
	If in addition $\alpha$ and $\beta$ satisfy the Cauchy-Riemann-equations (\ref{CR}), then $f$ is called right slice hyperholomorphic or right slice monogenic.
\end{definition}
\begin{definition}
	Let $U$ be an axially symmetric open set $\mathbb{R}^{n+1}$.
	\begin{enumerate}
		\item[(I)] We denote the set of left
		slice hyperholomorphic functions (or left slice monogenic functions) on $U$ by $\mathcal{SH}_L(U)$ and
		right slice hyperholomorphic functions (or right slice monogenic functions) on $U$ will be denoted by $\mathcal{SH}_R(U)$.\\
		\item[(II)] A slice hyperholomorphic function \eqref{form} such that $\alpha$ and $\beta$ are real-valued functions is called intrinsic slice hyperholomorphic or intrinsic slice monogenic functions and will be denoted by $\mathcal{N}(U)$.
	\end{enumerate}
\end{definition}
Function of the above type are systematically used in operator theory in the quaternionic and in the Clifford algebra setting,
see \cite{6COFBook,6ACSBOOK,FJBOOK,CGKBOOK,6JONAME}, and they are studied in the case of real alternative algebras in \cite{GP}.
\begin{remark}
	For quaternionic functions
	$f:U\subseteq \mathbb{H} \to \mathbb{H}$ the above definition is naturally re-adapted.
\end{remark}
We recall, for $s,\, x\in \mathbb R^{n+1}$ with $x\notin [s]$, the definition of the pseudo Cauchy kernel is given by
\begin{align*}
	\PRes_{c,s}(x)^{-1}:=(s^2-2\Re(x)s+|x|^2)^{-1},
\end{align*}
and the left slice hyperholomorphic Cauchy kernel  is defined as
\begin{equation}\label{Cauchyker}
	S^{-1}_L(s,x):=(s-\overline x) \PRes_{c,s}(x)^{-1},
\end{equation}
\begin{definition}\label{DIRAC}
	The Dirac operator $D$ and of its conjugate $\overline{D}$, are defined as
	\begin{equation}\label{DIRACeBARDNN}
		D= \frac{\partial}{\partial x_0}+ \sum_{i=1}^{n} e_i \frac{\partial}{\partial x_i}, \quad {\rm and}  \quad
		\overline{D}= \frac{\partial}{\partial x_0}- \sum_{i=1}^{n} e_i \frac{\partial}{\partial x_i}.
	\end{equation}
\end{definition}
\begin{definition}[Axially monogenic function]\label{AXIALM}
	Let us assume that $U \subseteq \mathbb{R}^{n+1}$ is an axially symmetric open set.
	A left (resp. right) axial function $f:U \to \mathbb{R}^{n+1}$ (see (\ref{form})) of class $ \mathcal{C}^1$ that satisfy the compatibility conditions (see \eqref{EO}) is said to be left (resp. right) axially monogenic if it is monogenic, i.e.
	\begin{align*}\label{axx}
		Df(x)=0 \qquad (\hbox{resp.} \quad 	f(x)D=0 ).
	\end{align*}
	The class of left (resp. right) axially monogenic is denoted by $ \mathcal{AM}_L(U)$ (resp. $\mathcal{AM}_R(U)$).
\end{definition}
On the sets of functions that admit an integral representation
deduced by the Cauchy formula of slice hyperholomorphic functions it is possible to define the related functional calculi whose resolvent operators are defined via the invertibility of the commutative version of the pseudo $S$-resolvent operator.
Precisely,  let $V$ be a real Banach space over $\rr$.
We denote by $\mathcal{B}(V)$ the space of all bounded $\mathbb{R}$-linear operators
and by $\mathcal{B}(V_n)$ the space of all bounded $\mathbb{R}_n$-linear operators
on $V_n=\mathbb{R}_n\otimes V$.
For the fine structures of the spectral theories on the $S$-spectrum
in $\mathcal{B}(V_n)$ we will consider bounded paravector operators
$T=e_0T_0+e_1T_1+...+e_n T_n$, with commuting components $T_\ell\in\mathcal{B}(V)$ for $\ell=0,1,\ldots ,n$
and we denote this set by $\mathcal{BC}^{\small 0,1}(V_n)$.
In this case the most appropriate definition of the $S$-spectrum is its commutative version (also called $F$-spectrum), i.e.,
\begin{align*}
	\sigma_S(T)=\{ s\in \mathbb{R}^{n+1}\ \ |\ \ s^2\mathcal{I}-(T+\overline{T})s +T\overline{T}\ \ \
	{\rm is\ not\  invertible\ in \ }\mathcal{B}(V_n)\}
\end{align*}
where the operator $\overline{T}$ is defined by
$\overline{T}=T_0-T_1e_1 - \dots  - T_n e_n.$
\medskip
Let $T\in \mathcal {BC}^{0,1}(V_n)$ and
recall that the $S$-resolvent set is defined as
$\rho_S(T):=\mathbb{R}^{n+1}\setminus \sigma_S(T),$
the commutative pseudo $S$-resolvent operator is given by
\begin{align*}
	\mathcal Q_{c,s}(T)^{-1}:=(s^2\mathcal I-s(T+\overline T)+T\overline T)^{-1} \quad  {\rm for} \quad   s\in \rho_S(T)
\end{align*}
while the left $SC$-resolvent operator is:
\begin{align*}
	S^{-1}_{L}(s,T):=(s\mathcal I-\overline T)(s^2\mathcal I-s(T+\overline T)+T\overline T)^{-1} \quad  {\rm for} \quad   s\in \rho_S(T).
\end{align*}


\section{The Fueter-Sce mapping theorem and spectral theories}\label{SEC2}

Starting from holomorphic functions, R. Fueter in 1934, see \cite{Fueter}, introduced a method to generate Cauchy-Fueter regular functions (i.e functions in the kernel of the operator $D$ with $n=3$). Over twenty years later, in 1957, M. Sce \cite{Sce} extended this result in a groundbreaking manner that contains Clifford algebras. For an English translation of his works in hypercomplex analysis, along with commentaries, refer to the recent book \cite{bookSCE}.
\begin{theorem}[Sce, 1957]
	Let $f$ be a holomorphic function in an open set of the upper half complex plane and let
	$f(x+iy)=u(x,y)+iv(x,y)$, for $x,y\in \mathbb{R}$,
	where $u$ and $v$ are real differentiable functions with values in $\mathbb{R}$.
	Consider the Euclidean space $\rr^{n+1}$ whose points $(x_0,....,x_n)$ are identified with
	the paravector $x_0+\underline{x}$ in the Clifford algebra $\rr_n$ and where  $\underline{x}={e_1}x_1+\ldots +{e_n}x_n$.
	Then
	the function
	$\breve{f}:=\textcolor{black}{\Delta_{n+1}^{\frac{n-1}{2}}} \Big(\textcolor{black}{u(x_0,|\underline{x}|)+\frac{\underline{x}}{|\underline{x}|}v(x_0,|\underline{x}|)}\Big)$
	is monogenic, i.e. it is in the kernel of the Dirac operator $D$.
\end{theorem}
\begin{remark}
	In the context of Clifford algebras, the operator $T_{FS2}$ is given by $ \Delta_{n+1}^{\frac{n-1}{2}}$, where $\Delta_{n+1}$ denotes the Laplace operator in $n+1$ variables. When $n$ is odd, $T_{F2}$ acts as a pointwise differential operator see \cite{Sce} and we will always consider this case. Conversely, if $n$ is even, the operator involves fractional powers of the Laplace operator, this was studied by T. Qian in \cite{TaoQian1}.
\end{remark}
The Fueter-Sce theorem is essential to understand the connections between hyperholomorphicity and the corresponding spectral theories that are based on hyperholomorphic functions.
Precisely, let $\mathcal{O}(\Pi)$ denote the set of holomorphic functions on $\Pi \subseteq \mathbb{C}$, and let $\Omega_\Pi \subseteq \mathbb{R}^{n+1}$ be the set induced by $\Pi$ and defined as $\Omega_\Pi:=\{x+Jy \, :\, (x,y)\in \Pi, J\in \mathbb{S}\}$. The first Fueter-Sce map $T_{FS1}$ applied to $\mathcal{O}(\Pi)$ produces the set ${\mathcal{SH}(\Omega_\Pi)}$ of slice monogenic functions on $\Omega_\Pi$, which are intrinsic. The second Fueter-Sce map $T_{FS2}$ applied to ${\mathcal{SH}(\Omega_\Pi)}$ generates axially monogenic functions $\mathcal{AM}(\Omega_\Pi)$. The extension procedure is illustrated as follows:
\begin{equation*}
	\begin{CD}
		\textcolor{black}{\mathcal{O}(\Pi)} @>T_{FS1}>> \textcolor{black}{\mathcal{SH}(\Omega_\Pi)} @> T_{FS2} = \Delta^{(n-1)/2}_{n+1} >> \textcolor{black}{\mathcal{AM}(\Omega_\Pi)}.
	\end{CD}
\end{equation*}
The Fueter-Sce mapping theorem induces two spectral theories based on the two classes of hyperholomorphic functions it generates. Using the Cauchy formula for slice hyperholomorphic functions defines the $S$-functional calculus based on the $S$-spectrum, while the Cauchy formula for monogenic functions leads to the monogenic functional calculus based on the monogenic spectrum. The latter calculus, introduced by A. McIntosh and collaborators \cite{JM}, is used to define functions of noncommuting operators on Banach spaces and has numerous applications, as discussed in \cite{J,TAOBOOK}.
This construction can also be visualized in the following diagram:
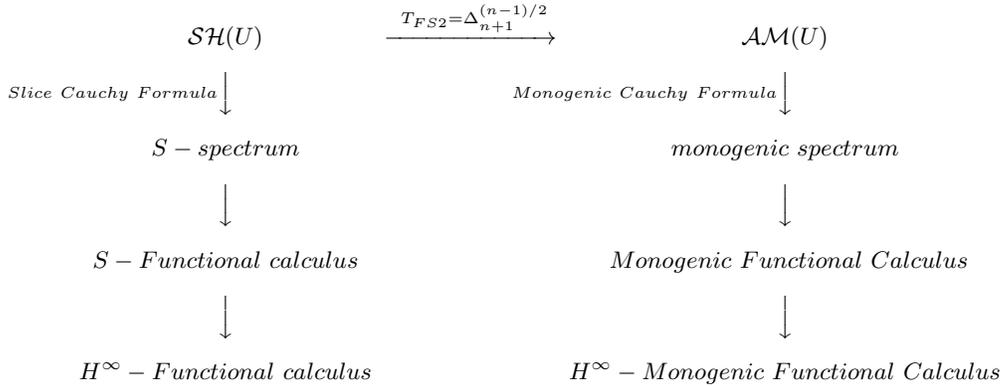
\begin{equation*}
	{\footnotesize
		\begin{CD}
			\textcolor{black}{\mathcal{SH}(U)} @> \ T_{FS2}=\Delta^{(n-1)/2}_{n+1} >>  \textcolor{black}{\mathcal{AM}(U)} \\   @V \textcolor{black}{Slice\ Cauchy \ Formula} VV
			@V \textcolor{black}{Monogenic \  Cauchy \ Formula} VV
			\\
			\textcolor{black}{S-spectrum}  @. \textcolor{black}{monogenic \ spectrum}
			\\
			@V VV    @V VV
			\\
			\textcolor{black}{S-Functional \ calculus} @. \ \textcolor{black}{Monogenic\ Functional\ Calculus}
			\\
			@V VV   @V VV
			\\
			\textcolor{black}{H^\infty-Functional \ calculus} @.   \textcolor{black}{H^\infty -Monogenic\ Functional\ Calculus}
		\end{CD}
	}
\end{equation*}
\begin{remark}
	The spectral theorem is based on the notion of the $S$-spectrum, see \cite{6SpecThm1} and \cite{CLIFST}.
\end{remark}
\begin{remark}
	In our diagrams we consider just the map $T_{FS2}$ from $\mathcal{SH}(\Omega_\Pi)$ to $\mathcal{AM}(\Omega_\Pi)$ and for this reason
	we will consider axially symmetric open sets, that will be denoted $U$, on which slice monogenic functions are defined so we will abandon the notation $\Omega_\Pi$.
\end{remark}


\section{What are the fine structures on the $S$-spectrum}

For the various functional calculi of the fine structure we mention the following
key references. The commutative version of the $S$-functional calculus is investigated in \cite{CS}, while the $F$-functional calculus is extensively discussed in \cite{CDS, CDS1, CGROYAL, CSS10}. Axially harmonic and harmonic functional calculi are elaborated in \cite{CDPS}. Polyanalytic functions and their associated functional calculus can be found in \cite{DP1,DP2}.
For unbounded operators on the quaternionic fine structure, see \cite{BANACHJ}. The $H^\infty$-functional calculus for fine structures is covered in \cite{CPSMILAN,DPSSPECTRAL}. Finally, the Clifford case in five dimensions is described in \cite{CDQS5dim}.

The Fueter-Sce mapping theorem (for odd \( n \)) offers an alternative approach to define the monogenic functional calculus. The core concept is applying the Fueter-Sce operator \( T_{FS2} \) to the slice hyperholomorphic Cauchy kernel and expressing the theorem in an integral form. This integral representation allows for the definition of the \( F \)-functional calculus, a monogenic functional calculus based on the \( S \)-spectrum. The process can be illustrated as follows:
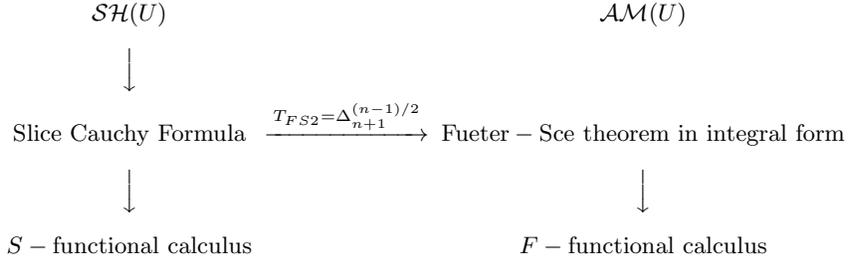
\begin{equation*}
	{\footnotesize
		\begin{CD}
			{\mathcal{SH}(U)} @.  {\mathcal{AM}(U)} \\   @V  VV
			@.
			\\
			{{\rm  Slice\ Cauchy \ Formula}}  @> T_{FS2}=\Delta^{(n-1)/2}_{n+1}>> {{\rm Fueter-Sce\ theorem \ in \ integral\ form}}
			\\
			@V VV    @V VV
			\\
			{S-{\rm  functional \ calculus}} @. F-{{\rm functional \ calculus}}
		\end{CD}
	}
\end{equation*}
Note that the diagram lacks an arrow from the axially monogenic function space \(\mathcal{AM}(U)\) because the \( F \)-functional calculus is derived from the slice hyperholomorphic Cauchy formula.
Specifically, for odd $n$ in the Clifford algebra $\mathbb{R}_n$, one of the most significant factorizations leads to the so called Dirac fine structure.
This structure involves a suitable alternating sequence of the operators $D$ and $\overline{D}$
repeated $(n-1)/2$ times.
As an example of a factorization of the operator $T_{FS2}$ we have:
\begin{align*}
	T_{FS2}=\Delta_{n+1}^{h}= {D} {\overline{D}}\cdots {D}{\overline{D}} .
\end{align*}
where $h:=(n-1)/2$ is the Sce exponent.
Since $D$ and $\overline{D}$ commute with each other,
various configurations of their order can be considered, as long as their product results in $\Delta_{n+1}^{h}$.
As a further example also the configuration
\begin{align*}
	\Delta_{n+1}^{h}= {D} \cdots {D} {\overline{D}}  \cdots{\overline{D}}
\end{align*}
is possible
and if we allow in certain positions $D\overline{D}$ or $ \overline{D}D$
to be the Laplacian $\Delta_{n+1}$, since
$\Delta_{n+1}=D\overline{D}=\overline{D}D$
several different fine structures emerge from the operators
$D$, $\overline{D}$, $\Delta_{n+1}$ and their powers, depending on the dimension of the Clifford algebra.

We are now in the position to define the fine structure of spectral theories on the $S$-spectrum, utilizing the above observation.
\begin{definition}
	We will call
	{\rm fine structure of the  spectral theory on the $S$-spectrum}
	\begin{itemize}
		\item
		{\em the set of the function spaces} and
		\item
		{\em the associated functional calculi}
	\end{itemize}
	induced by all possible factorizations of the operator $T_{FS2}$ in the Fueter-Sce extension theorem.
\end{definition}


\subsection{A unified way to see some function spaces}\label{UNIFFWAY}

The fine structure of the spectral theory on the $S$-spectrum
generates, in a unified way, several classes of functions,
some of which have already been studied in the literature.
We can summarise them by defining  polyanalytic holomorphic Cliffordian functions of order $(k, \ell)$ which are defined as follows.
Because of the factorization
${D}{\overline{D}}={\overline{D}}{D}=\Delta_{n+1}$
of the Laplace operator
we define some classes of functions that are strictly related to the Fueter-Sce theorem.
We will see in the next sections that these function spaces are of crucial importance for our theory.
\begin{definition}[holomorphic Cliffordian of order $k$]
	\label{hc}
	Let $U\subset \mathbb{R}^{n+1}$ be an open set and let  $0 \leq k \leq (n-1)/2$ be an integer. A function
	$f:U \subset \mathbb{R}^{n+1} \to \mathbb{R}_n$ of class $\mathcal{C}^{2k+1}(U)$ is said to be (left) holomorphic Cliffordian of order $k$ if
	$ \Delta_{n+1}^k D  f(x)=0,$ $\forall x\in U.$
\end{definition}
\begin{remark} We observe that
	for $k:= (n-1)/2$ in Definition \ref{hc} we get the class of functions studied in \cite{LR}.
	Moreover, every holomorphic Cliffordian function of order $k$ is holomorphic Cliffordian of order $k+1$. If $k=0$ in Definition \ref{hc} we get the set of (left) monogenic functions.
\end{remark}
\begin{definition}[anti-holomorphic Cliffordian of order $k$]
	\label{anti1}
	Let $U\subset \mathbb{R}^{n+1}$ be an open set and let  $0 \leq k \leq (n-1)/2$ be an integer. A function $f:U \subset \mathbb{R}^{n+1} \to \mathbb{R}_n$ of class $\mathcal{C}^{2k+1}(U)$ is said to be (left) anti-holomorphic Cliffordian of order $k$ if 
	$ \Delta_{n+1}^k \overline{D}  f(x)=0$, $\forall x\in U.$
\end{definition}
\begin{definition}[polyharmonic of degree $k$]
	\label{harmo}
	Let $k \geq 1$. A function $f:U \subset \mathbb{R}^{n+1} \to \mathbb{R}_n$ of class $\mathcal{C}^{2k}(U)$ is called polyharmonic of degree $k$ in the open set $U \subset \mathbb{R}^{n+1}$ if for any $x \in U$ we have
	$ \Delta_{n+1}^k f(x)=0.$
\end{definition}
For $k=1$ the function is called harmonic and for $k=2$ the function is called bi-harmonic.
\begin{definition}[polyanalytic of order $m$, see \cite{B1976}]
	\label{poly}
	Let $ m \geq 1$. Let $ U \subset \mathbb{R}^{n+1}$ be an open set and let $f:U \to \mathbb{R}_n$ be a function of class $ \mathcal{C}^{m}(U)$. We say that $f$ is (left) polyanalytic of order $m$ on $U$ if for any $x \in U$ we have
	$ D^{m}f(x)=0.$
\end{definition}
The above spaces, apart from anti-holomorphic Cliffordian of order $k$,
have already been considered in the literature and can be summarized in the following definition.
\begin{definition}[polyanalytic holomorphic Cliffordian of order $(k, \ell)$]
	\label{polycl}
	Let \linebreak $U\subset \mathbb{R}^{n+1}$ be an open set and let  $0 \leq k \leq (n-1)/2$ be an integer.
	A function $f:U \subset \mathbb{R}^{n+1} \to \mathbb{R}_n$ of class $\mathcal{C}^{2k+ \ell}(U)$ is said to be (left) polyanalytic holomorphic Cliffordian of order $(k, \ell)$ if for any $x \in U$ we have
	$ \Delta_{n+1}^k D^{\ell}  f(x)=0$. We denote the set of these functions as $\mathcal{PCH}_{(k,\ell)}(U)$.
\end{definition}
\begin{remark}
	Similarly it is possible to define the  sets of right functions for the holomorphic Cliffordian of order $k$, anti-holomorphic Cliffordian of order $k$ and polyanalytic of order $m$.
\end{remark}
\begin{remark}
	If in Definition \ref{polycl} we set $ \ell=1$ we get Definition \ref{hc}, if $\ell=0$ we obtain Definition \ref{harmo} and if we consider $k=0$ we obtain Definition \ref{poly}.
\end{remark}
\medskip \medskip
Clearly, as the dimension of the Clifford algebra increases
there are more possibilities for the
functions spaces that lie between the set of slice hyperholomorphic functions and axially monogenic functions.
Moreover, it is very important to point out that these function spaces
appear in different contexts
in the literature and they seem to be unrelated, but the fine structures on the $S$-spectrum give them a unified setting.


\section{The quaternionic setting}

The quaternionic fine structures of Dirac type is now based on the two different ways we can factorize the Laplacian, in dimension four,
$\Delta_4= D\overline{ D}=\overline{ D} D,$
using the Cauchy-Fueter operator (also called Dirac-operator) $D$ and its conjugate
$\overline{ D}$ are as in Definition \ref{DIRAC} with $n=4$.
In the quaternionic setting, using the factorization of the second map in the Fueter mapping theorem, that is the Laplacian in four dimensions, we get two fine structures.
In \cite{CDPS} we studied the Dirac fine structure when $n=3$, i.e., the quaternionic case.
Depending on whether $ D$ or $\overline{D}$ is applied first on some function $f\in\mathcal{SH}(U)$, we get the following four function spaces:
\begin{remark}[The function spaces of the quaternionic fine structure]\label{PRECISEDEF}
	Let $\mathcal{SH}(U)$ be as in Definition \ref{sh}. Then, we define
	\begin{itemize}
		\item
		$\mathcal{AH}(U)=\{ Df | f\in\mathcal{SH}(U)\}$,
		{\rm (axially harmonic functions)}
		\item
		$\mathcal{AP}_2(U)=\{\overline{ D}f | f\in\mathcal{SH}(U)\}$,
		{\rm(polyanalytic functions of order 2)}
		\item
		$\mathcal{AM}(U)=\{\Delta_4 f | f\in\mathcal{SH}(U)\},$
		{\rm (axially monogenic functions).}
	\end{itemize}
\end{remark}
\begin{remark}\label{COMMENT CK}
	We note that, with some slight abuse of notation, in Remark \ref{PRECISEDEF} we will keep
	the names of the function spaces defined in the previous Subsection \ref{UNIFFWAY}, that
	only require $\mathcal{C}^k$ regularity for suitable $k \in \mathbb{N}.$
	Within the context of the Fueter-Sce theorem
	the same consideration have to be done for the Clifford algebra setting.
\end{remark}
In the quaternionic setting we have studied the fine structure associated with the factorization:
\begin{equation}
	\label{fine1} \mathcal{SH}(U)\overset{D}{\longrightarrow} \mathcal{AH}(U)\overset{\overline{D}}{\longrightarrow}\mathcal{AM}(U),
\end{equation}
where $\mathcal{AH}(U)$ is the set of axially harmonic functions
and their integral representation give rise to the harmonic functional calculus on the $S$-spectrum.
This structure also allows to obtain a product formula for the $F$-functional calculus, see \cite[Thm. 9.3]{CDPS}.
However, since $\Delta_4={D} {\overline{D}}= {\overline{D}}{D}$, we can interchange the order of the operators $ D$ and $ \overline{D}$ in \eqref{fine1}. This gives rise to the factorization:
\begin{equation} \mathcal{SH}(U)\overset{\overline{D}}{\longrightarrow} \mathcal{AP}_2(U)\overset{D}{\longrightarrow}\mathcal{AM}(U),
\end{equation}
where $\mathcal{AP}_2(U)$ is a space of polyanalytic functions. This structure is investigated in \cite{DP1} and \cite{DP2} together with its functional calculus.
The harmonic fine structure of the quaternionic spectral
theory on the $S$-spectrum is illustrated in the following diagram
{\footnotesize
	\textcolor{black}
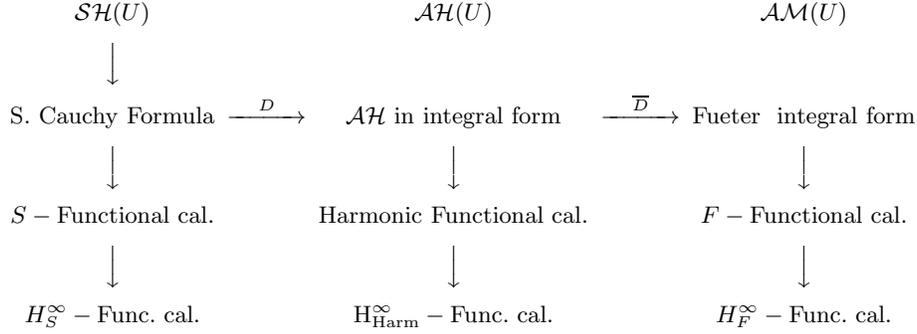
{\begin{equation*}
			\begin{CD}
				{\mathcal{SH}(U)} @. {\mathcal{AH}(U)}  @.  {\mathcal{AM}(U)} \\   @V  VV
				@.
				\\
				{{\rm  S.\ Cauchy \ Formula}}  @> D >> {\mathcal{AH} {\rm \ in \  integral\  form}}@> \overline{D} >> {{\rm Fueter\ \ integral\  form}}
				\\
				@V VV    @V VV  @V VV
				\\
				S-{{\rm Functional \ cal.}} @. {{\rm Harmonic \ Functional \ cal.}}@. F-{{\rm Functional \ cal.}}\\
				@V VV    @V VV  @V VV
				\\
				H_S^\infty -{{\rm Func.\ cal.}} @. {{\rm H_{Harm}^\infty - Func. \ cal.}}@. H_F^\infty -{{\rm Func. \ cal.}}
			\end{CD}
		\end{equation*}
}}
where the central part of the diagram contains the harmonic functions and their functional calculi.
The polyanalytic fine structure of the quaternionic spectral theory on the $S$-spectrum is illustrated in the following diagram
{\footnotesize
	\textcolor{black}
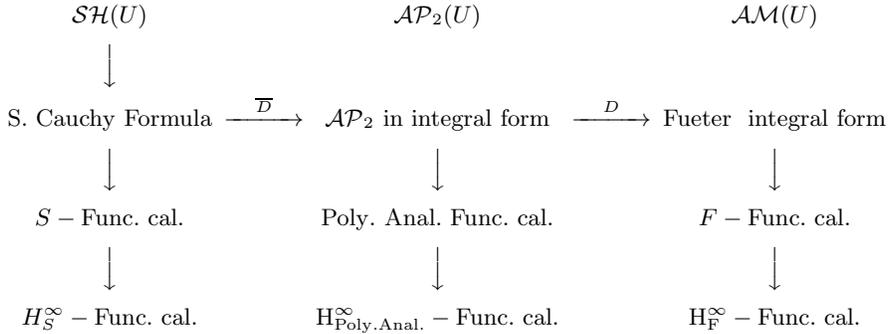
{\begin{equation*}
			\begin{CD}
				{\mathcal{SH}(U)} @. {\mathcal{AP}_2(U)}  @.  {\mathcal{AM}(U)} \\   @V  VV
				@.
				\\
				{{\rm  S.\ Cauchy \ Formula}}  @> \overline{D} >> {\mathcal{AP}_2 {\rm \ in \  integral\  form}}@> {D} >> {{\rm Fueter\ \ integral\  form}}
				\\
				@V VV    @V VV  @V VV
				\\
				S-{{\rm Func. \ cal.}} @. {{\rm Poly.\ Anal. \ Func. \ cal.}}@. F-{{\rm Func.\ cal.}}\\
				@V VV    @V VV  @V VV
				\\
				H_S^\infty -{{\rm Func. \ cal.}} @.  {{\rm H_{Poly. Anal.}^\infty- Func.\ cal.}}@.  {{\rm H_F^\infty-Func. \ cal.}}
			\end{CD}
		\end{equation*}
}}
\newline
where  the central part of the diagram  contains the polyanalytic functions and their functional calculi.


\section{The five dimensional case}

In a Clifford algebra with $n$ imaginary units, there are numerous possibilities for the function spaces, and the study of fine structures is advancing rapidly.
The case of dimension five already shows that most of the function spaces
discussed in the previous section appear.
Below, we provide several examples to illustrate this theory, which involves the spaces explicitly listed below, keeping in mind Remark \ref{COMMENT CK}:
\begin{itemize}
	\item
	$ \mathcal{ABH}(U)$ the axially bi-harmonic functions (\textcolor{black}{$\Delta_6 ^2f(x)=0$}),
	\item
	$ \mathcal{ACH}_1(U)$ the axially Cliffordian holomorphic functions of order 1,
	(which is a short cut for order $(1,1)$), (\textcolor{black}{$ \Delta_6 D f(x)=0$}),
	\item
	$ \mathcal{AH}(U)$ the axially harmonic functions, (\textcolor{black}{$\Delta_6 f(x)=0$})
	\item
	$ \mathcal{AP}_2(U)$ the axially polyanalytic of order $2$, (\textcolor{black}{$D^2 f(x)=0$})
	\item
	$ \mathcal{ACH}_1(U)$ the axially anti Cliffordian of order $1$,
	(\textcolor{black}{$\overline{D} f(x)=0)$})
	\item
	$ \mathcal{ACP}_{(1,2)}$ the axially polyanalytic Cliffordian of order $(1,2)$,
	(\textcolor{black}{$ \Delta_6 D^2 f(x)=0$}),
	\item
	$ \mathcal{AP}_3(U)$ the axially polyanalytic of order $3$,
	(\textcolor{black}{$D^3 f(x)=0$}).
\end{itemize}
Rearranging the sequence of $ D$ and $ \overline{D}$ it is possible to obtain different fine structures, in which the above sets of functions are involved. Thus, we have the Dirac fine structure  $(D, \overline{D},\overline{D}, D )$
\begin{align*}
	\mathcal{SH}(U)\overset{D}{\longrightarrow} \mathcal{ABH}(U)\overset{\overline{D}}{\longrightarrow}	\mathcal{AHC}_1(U) \overset{\overline{D}}{\longrightarrow} \mathcal{AP}_2(U) \overset{D}{\longrightarrow} \mathcal{AM}(U),
\end{align*}
and the Dirac fine structure $(D, D, \overline{D}, \overline{D}  )$
\begin{align*}
	\mathcal{SH}(U)\overset{D}{\longrightarrow} \mathcal{ABH}(U)\overset{D}{\longrightarrow}\overline{\mathcal{AHC}_1(U)} \overset{\overline{D}}{\longrightarrow} \mathcal{AH}(U) \overset{\overline{D}}{\longrightarrow} \mathcal{AM}(U).
\end{align*}
All the previous Dirac fine structures are obtained by applying first the Dirac operator. Nevertheless, it is possible to apply the operator $ \overline{D}$ as  first operator. In this case other three Dirac fine structures arise.
We have the Dirac fine structure of the kind $(\overline{D}, D, \overline{D},D )$
\begin{align*}
	\mathcal{SH}(U)\overset{\overline{D}}{\longrightarrow} \mathcal{APC}_{(1,2)}(U)\overset{D}{\longrightarrow} \mathcal{AHC}_1(U) \overset{\overline{D}}{\longrightarrow} \mathcal{AP}_2(U) \overset{D}{\longrightarrow} \mathcal{AM}(U),
\end{align*}
the Dirac fine structure $(\overline{D}, D,D,  \overline{D} )$
\begin{align*}
	\mathcal{SH}(U)\overset{\overline{D}}{\longrightarrow} \mathcal{APC}_{(1,2)}(U)\overset{D}{\longrightarrow} \mathcal{AHC}_1(U) \overset{D}{\longrightarrow} \mathcal{AH}(U) \overset{\overline{D}}{\longrightarrow} \mathcal{AM}(U),
\end{align*}
and the Dirac fine structure $(\overline{D},\overline{D}, D,D)$
\begin{align*}
	\mathcal{SH}(U)\overset{\overline{D}}{\longrightarrow} \mathcal{APC}_{(1,2)}(U)\overset{\overline{D}}{\longrightarrow} \mathcal{AP}_3(U) \overset{D}{\longrightarrow} \mathcal{AH}(U) \overset{D}{\longrightarrow} \mathcal{AM}(U).
\end{align*}
In all the previous Dirac fine structures it is possible to combine the Dirac operator and its conjugate. In this way we get a fine structure which is "weaker" then the previous ones; in the sense that we are skipping some classes of functions.

For example the Laplace fine structure is of the kind $(\Delta, \Delta)$,
which is a coarser fine structure with respect to the Dirac one, and it can be represented by the following diagram
\begin{align*}
	\mathcal{SH}(U)\overset{\Delta_6}{\longrightarrow}  \mathcal{ACH}_1(U) \overset{\Delta_6}{\longrightarrow}  \mathcal{AM}(U).
\end{align*}
Another interesting coarser fine structure is the polyanalytic one, in which there appear only the polyanalytic functions of order three and two
\begin{align*}
	\mathcal{SH}(U)\overset{\overline{D}^2}{\longrightarrow}  \mathcal{AP}_3(U) \overset{D}{\longrightarrow}  \mathcal{AP}_2(U)\overset{D}{\longrightarrow} \mathcal{AM}(U).
\end{align*}
\medskip
An additional intriguing coarser fine structure is the harmonic one, in which only the harmonic $\mathcal{AH}(U)$ and bi-harmonic $\mathcal{ABH}(U)$ sets of functions appear:
\begin{align*}
	\mathcal{SH}(U)\overset{D}{\longrightarrow}  \mathcal{ABH}(U) \overset{\Delta_6}{\longrightarrow}  \mathcal{AH}(U)\overset{\overline{D}}{\longrightarrow} \mathcal{AM}(U).
\end{align*}

The integral  representations of the functions in $\mathcal{AH}(U)$ and $\mathcal{ABH}(U)$ is obtained
by applying the operators in the above sequence to the
Cauchy kernels of slice hyperholomorphic functions.
From these integral formulas
we define the related functional calculi on the $S$-spectrum that can be visualized by the following diagram:
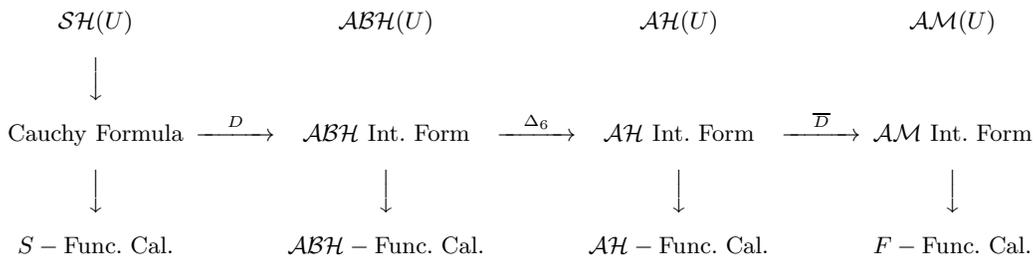
\begin{equation*}
	{\footnotesize
		\begin{CD}
			{\mathcal{SH}(U)}  @. \mathcal{ABH}(U)  @. \mathcal{AH}(U)  @.{\mathcal{AM}(U)} \\   @V  VV
			@.
			\\
			{{\rm  Cauchy\ Formula}}  @> D>> {{\rm \mathcal{ABH}\ \mathrm{Int.\ Form}}} @> \Delta_6>> {{\rm
					\mathcal{AH}\ \mathrm{Int.\ Form}}}@> \overline{D}>> {{\rm \mathcal{AM}\ \mathrm{Int.\ Form}  }}
			\\
			@V VV    @V VV  @V VV    @V VV
			\\
			{S-{\rm  \mathrm{Func.\ Cal.} }} @. {{\rm \mathcal{ABH}-\mathrm{Func.\ Cal.} }}@. {{\rm \mathcal{AH} -\mathrm{Func.\ Cal.}}}@. {F-{\rm \mathrm{Func.\ Cal.}}}
		\end{CD}
	}
\end{equation*}
where the functional calculi of this fine structure are the
bi-harmonic $\mathcal{ABH}$-functional calculus and the harmonic $\mathcal{AH}$-functional calculus,
both based on the $S$-spectrum.

We observe that it is not possible to have coarser fine structure in the quaternionic case. This is due to the fact that we are dealing with the Laplacian of power $1$.
\medskip\medskip

We denote by  $\mathcal{SH}_L(\sigma_S(T))$
the set of functions $f\in \mathcal{SH}_L(U)$ where $U$ is a bounded slice Cauchy domain with $\sigma_S(T)\subset U$ and $\overline{U}\subset\operatorname{dom}(f)$.

Considering the Cauchy kernel (\ref{Cauchyker}) we define the kernels 
$S^{-1}_{\Delta_6^{1-\ell}D,L}(s,x):=\Delta_6^{1-\ell}D S^{-1}_{L}(s,x)$.
Then, the resolvent operators $S^{-1}_{\Delta_6^{1-\ell} D,L}(s,T)$ are defined
replacing $x$ by $T$ in $\Delta_6^{1-\ell} DS^{-1}_{L}(s,x)$.
\begin{definition}
	Let $n=5$, assume $T \in \mathcal{BC}^{0,1}(V_5)$ and set $ds_I=ds(-I)$ for $I \in \mathbb{S}$.
	For $f$ in $\mathcal{SH}_L(\sigma_S(T))$
	we define the $\ell+1$-harmonic functional calculus for $\ell=0,1$, for every function $\tilde{f}_\ell=\Delta_6^{1-\ell}Df$ with $f \in \mathcal{SH}_L(\sigma_S(T))$, as
	\begin{equation}
		\tilde{f}_\ell(T):=\frac{1}{2\pi} \int_{\partial(U \cap \mathbb{C}_I)} S^{-1}_{\Delta_6^{1-\ell}D,L}(s,T)ds_I f(s),
	\end{equation}
	where
	\begin{itemize}
		\item the left $D$-resolvent operator $S^{-1}_{D,L}(s,T)$ is defined as
		\begin{equation}
			S^{-1}_{D,L}(s,T):=-4 \mathcal{Q}_{c,s}(T)^{-1} \ \ {\rm for}\ \ \ s\in \rho_S(T),
		\end{equation}
		\item the left $\Delta_6 D$-resolvent operator $S^{-1}_{\Delta_6 D,L}(s,T)$ is defined as
		\begin{equation}
			S^{-1}_{\Delta_6 D,L}(s,T):=16 \mathcal{Q}_{c,s}(T)^{-2}, \ \ {\rm for}\ \ \ s\in \rho_S(T).
		\end{equation}
	\end{itemize}
\end{definition}

\medskip We point out that similar formulas hold for the right case.
\newline
\newline
\newline
\textbf{Acknowledgement}
	Peter Schlosser was funded by the Austrian Science Fund (FWF) under Grant No. J 4685-N and by the European Union--NextGenerationEU.

\end{document}